\documentclass[12pt]{article}
\usepackage{amsmath}
\usepackage{amssymb}
\usepackage{amsthm}
\usepackage{epsfig}
\usepackage{graphicx}
\usepackage{enumerate}
\usepackage{cite}
\usepackage{color}

\theoremstyle{remark}

\numberwithin{equation}{section}

\begin{document}

\title{Tables of sizes of random complete arcs in the plane $PG(2,q)$}
\date{}
\author{Daniele Bartoli \\
{\footnotesize Dipartimento di Matematica e Informatica,
Universit\`{a}
degli Studi di Perugia, }\\
{\footnotesize Via Vanvitelli~1, Perugia, 06123, Italy.
E-mail: daniele.bartoli@dmi.unipg.it}\\
 \and Alexander A.
Davydov \\
{\footnotesize Institute for Information Transmission Problems
(Kharkevich
institute), Russian Academy of Sciences}\\
{\footnotesize Bol'shoi Karetnyi per. 19, GSP-4, Moscow,
127994, Russia. E-mail: adav@iitp.ru}
\and Giorgio Faina, Stefano Marcugini and Fernanda Pambianco \\
{\footnotesize Dipartimento di Matematica e Informatica,
Universit\`{a}
degli Studi di Perugia, }\\
{\footnotesize Via Vanvitelli~1, Perugia, 06123, Italy. E-mail:
\{faina,gino,fernanda\}@dmi.unipg.it} }\maketitle

\begin{abstract} Tables of sizes of random complete arcs in the plane $PG(2,q)$ are given. The sizes
 are  close to the smallest known
sizes of complete arcs in $PG(2,q)$, in particular, to ones
constructed by Algorithm FOP (fixed order of points). The
random arcs are obtained in the region $\{3\leq q\leq 46337,~q
\mbox{ prime}\}$.
\end{abstract}



\section{Introduction}

Let $PG(2,q)$ be the projective plane over the Galois field
$F_{q}$. An $n$-arc is a set of $n$ points no three of which
are collinear. An $n$-arc is called complete if it is not
contained in an $(n+1)$-arc of $PG(2,q)$. For an introduction
to projective geometries over finite fields see \cite
{HirsBook,SegreLeGeom,SegreIntrodGalGeom}.

In \cite{tallini,HirsStor-2001} the close relationship among
the theory of $n$-arcs, coding theory, and mathematical
statistics is presented. In particular, a complete arc in a
plane $PG(2,q),$ the points of which are treated as
3-dimensional $q$-ary columns, defines a parity check matrix of
a $q$-ary linear code with codimension 3, Hamming distance 4,
and covering radius 2. Arcs can be interpreted as linear
maximum distance separable (MDS) codes \cite[Sec. 7]{szoT93},
\cite{thaJ92d} and they are related to optimal coverings arrays
\cite{Hartman-Haskin} and to superregular matrices \cite
{Keri}.

One of the main problems in the study of projective planes,
which is also of interest in coding theory, is finding of the
spectrum of possible sizes of complete arcs. In particular, the
value of $t_{2}(2,q)$, the smallest size of a complete arc in
$PG(2,q),$ is interesting. Finding estimates of the minimum
size $t_{2}(2,q)$ is a hard open problem.

This work is devoted to \emph{random complete arcs} in
$PG(2,q)$ and to the comparison of their sizes with \emph{upper
bounds} on $t_{2}(2,q)$.

Surveys of results on the sizes of plane complete arcs, methods
of their construction, and comprehension of the relating
properties can be found in
\cite{BDFMP-JG2013,BDFKMP-Bulg2013,BDFMP2013,BDFMP-ArXivFOP,BDFMP_JGsubmit,BDFMP2012,DFMP-JG2005,DFMP-JG2009,FP,HirsSurvey83,HirsBook,HirsStor-2001,%
KV,LombRad,pelG77,pelG92,pel93,SegreLeGeom,segB62b,SegreIntrodGalGeom,SZ,szoT87a,szoT89survey,szoT93,Szonyi1997surveyCurves,PBFM2013,BFMP2013b}.
Some problems connected with small complete plane arcs are
considered in
\cite{abaV83,Ball-SmallArcs,BDFKMP-Bulg2013,BDFKMP-ArXiv,BDFMP-ArXivFOP,BDFMP_JGsubmit,DGMP-JCD,DGMP-AMC,FPDM,GacsSzonyi,MaPa-private,GiulUghi,%
HirsSurvey83,HirsBook,HirsStor-2001,KV,korG83a,Ost,PDBGM2013,Polv,Voloch90}.

The exact values of $t_{2}(2,q)$ are known only for $q\leq 32$;
see \cite
{tesi13,faina24,faina15,HirsBook,6,faina40,MMP-q29,art25} and
 work \cite{BFMP2013} where the equalities $
t_{2}(2,31)=t_{2}(2,32)=14$ are established.

Let $t(\mathcal{P}_{q})$ be the size of the smallest complete
arc in any (not necessarily Galois) projective plane
$\mathcal{P}_{q}$ of order $q$. In \cite{KV}, for
\emph{sufficiently large} $q$, the following result is proven
by \emph{probabilistic methods} (we give it in the form of
\cite[ Table 2.6]{HirsStor-2001} taking into account that all
logarithms in \cite{KV} have natural base, see
\cite[p.\thinspace 10]{KV}):
\begin{equation}
t(\mathcal{P}_{q})\leq D\sqrt{q}\ln^{C}q,\text{ }C\leq 300,
\label{eq1_KimVu_c=300}
\end{equation}
where $C$ and $D$ are constants independent of $q$ (so-called
universal or absolute constants).
 The
authors of \cite{KV} conjecture that the constant can be
reduced to $C=10$. A survey and an analysis of random
constructions for geometrical objects can be found in \cite
{GacsSzonyi}; see also the references therein.

 Regarding complete arcs of sizes smaller $\frac{1}{2}q$
obtained by algebraic constructions, following
\cite[p.\thinspace209] {HirsStor-2001}, complete arcs in
$PG(2,q)$ have been constructed with sizes approximately
$\frac{1}{3}q$ (see \cite
{abaV83,BDFMP2012,korG83a,SZ,szoT87a,Voloch90}), $\frac{1}{4}q$
(see \cite {BDFMP2012,korG83a,szoT89survey}), $2q^{0.9}$ (see
\cite{SZ} where such arcs are constructed for $q>7^{10}$). It
is noted in \cite[Sec. 8] {GacsSzonyi}, that the smallest size
of a complete arc in $PG(2,q)$ obtained via algebraic
constructions is $cq^{3/4}$ where $c$ is a universal constant;
see \cite[Sec. 3]{szoT89survey} and \cite[Th. 6.8]{szoT93}.

In \cite{BDFMP2012,BDFMP-JG2013},  for large ranges of $q$, the
form of the bound of \eqref{eq1_KimVu_c=300} is applied but the
value of the constant $C$ was essentially reduced to $C=0.75$
\cite{BDFMP2012} and to $C=0.73$ \cite{BDFMP-JG2013} whereas
$D<1$. In particular, the following results are obtained in
\cite{BDFMP2012,BDFMP-JG2013} using  randomized greedy
algorithms:
\begin{eqnarray}\label{eq1_C-constant}
t_{2}(2,q) &<& \sqrt{q}\ln ^{0.75}q\quad \mbox{ for
}~23\leq q\leq ~~5107  \mbox{ \cite{BDFMP2012}};\\
t_{2}(2,q) &<&\sqrt{q}\ln ^{0.73}q\quad \mbox{ for }109\leq
q\leq 13627\mbox{ \cite{BDFMP-JG2013}}.
\end{eqnarray}

In \cite{BDFKMP-ArXiv}, the smallest known sizes of complete
arcs in $PG(2,q)$ (up to November 2013) are collected for the
following huge region $H$:
\begin{align} &H=\{q: 173\le q\le49727,~
q\mbox{ power prime}\}\cup\{q:173\le q\leq
125003,~ q \mbox{ prime}\}\,\cup  \label{eq1_R} \\
& \{59 \mbox{ sporadic prime $q$'s in the interval }\lbrack
125101\ldots360007], \mbox{ see \cite[Table
7]{BDFKMP-ArXiv}}\}. \notag
\end{align}
The data collected in
\cite{BDFKMP-ArXiv,BDFMP2012,BDFMP-JG2013} provide the
following result.
\begin{align}
t_{2}(2,q)<\sqrt{q}\ln ^{0.7295}q\quad \text{for }109\leq q\leq 169 \text{
and }q\in H.
\end{align}

In the recent works of the authors, see
\cite{BDFMP_JGsubmit,BFMP2013b,BDMP_ACCT2012}, a new Algorithm
FOP (fixed order of points) constructing small complete arcs in
$PG(2,q)$ is proposed. Lexicographical and the Singer fixed
orders of points are investigated.  We denote
\begin{align}\label{eq1_primeLex}
&L=\{q:3\le q\le67993,~ q \mbox{ prime}\}\,\cup
 \{43 \mbox{ sporadic prime $q$'s in }[69997\ldots190027]\};\\
 &S=\{q:5\le q\le40009,~ q \mbox{ prime}\}.\label{eq1_primeSing}
\end{align}
Let  $t_{2}^{L}(2,q)$ be the size of complete arcs in $PG(2,q)$
obtained by Algorithm FOP with Lexicographical order of points.
Let  $t_{2}^{S}(2,q)$ be the size of complete arcs in $PG(2,q)$
obtained by Algorithm FOP with Singer order of points. Values
of $t_{2}^{L}(2,q)$ in the region $L$ and $t_{2}^{S}(2,q)$ in
the region $S$ are collected in \cite{BDFMP-ArXivFOP}.

The data collected in \cite{BDFMP-ArXivFOP} provide the
following \textbf{upper bounds} on $t_{2}(2,q)$:
\begin{align}
&t_{2}(2,q)<t_{2}^{L}(2,q)<1.83\sqrt{q\ln q}\mbox{~ if ~} q\in L;\label{eq1_LexBound} \\
&t_{2}(2,q)<t_{2}^{S}(2,q)<1.83\sqrt{q\ln q}\mbox{~ if ~} q\in S.\label{eq1_SingBound}
\end{align}

In \cite{Blokhuis,KV} it is noted that, in a preliminary report
in 1989, J.\ C.\ Fisher  obtained by computer search complete
arcs in many planes of small orders and conjectured that
average size of a complete arc is about
\begin{align}\label{eq1_aversize}\sqrt{ 3q\log
q}.
\end{align}

We denote
\begin{align}\label{eq1_regionR}R=\{3\leq q\leq 46337,~q \mbox{ prime}\}.
\end{align}

 \emph{In this work}, we collect the sizes
$t_{2}^{R}(2,q)$ of \emph{random complete arcs} in $PG(2,q)$ in
the region~$R$. The collected sizes are represented in Table 1
 and in Figure \ref{fig1}. For comparison, we also give
Figures \ref{fig2} and \ref{fig3} with the sizes
$t_{2}^{L}(2,q)$ and $t_{2}^{S}(2,q)$ of complete arcs obtained
by Algorithm FOP with Lexicographical and Singer orders of
points. Finally, we represent differences
$t_{2}^{L}(2,q)-t_{2}^{R}(2,q)$ and
$t_{2}^{S}(2,q)-t_{2}^{R}(2,q)$ in Figures \ref{fig4} and
\ref{fig5}.

The random arcs are obtained in this work with the help of a
random generator used in a C++ program under the system Linux.
A complete arc is constructed step-by-step in a random manner.
At every step a point of the plane is selected randomly: if the
point is not covered by bisecants of the arc, then it is added
to the arc; otherwise, another point is selected. The process
stops when a complete arc is obtained.

From Table 1 and Figure \ref{fig1} it follows that
\begin{align}
t_{2}^{R}(2,q)<1.83\sqrt{q\ln q}\mbox{~ if ~} q\in R.\label{eq1_SatisfBound}
\end{align}
So, \emph{the sizes of random arcs in the region $R$ satisfy
the upper bounds on $t_{2}(2,q)$ given in}
\eqref{eq1_LexBound}, \eqref{eq1_SingBound}. One can say also
that \emph{the conjecture \eqref{eq1_aversize} holds in the
region $R$, see Figure} \ref{fig1}.

\section{Table. Figures}

\hspace{0.5cm}The sizes
 $t^{R}_{2}=t^{R}_{2}(2,q)$
 of random complete arcs in planes
 $\mathrm{PG}(2,q),$ $3\leq q\leq 46337$, $q$ prime, are shown in Table 1, see pp.
 10--32.

\begin{figure}[tbp]
\epsfig{file=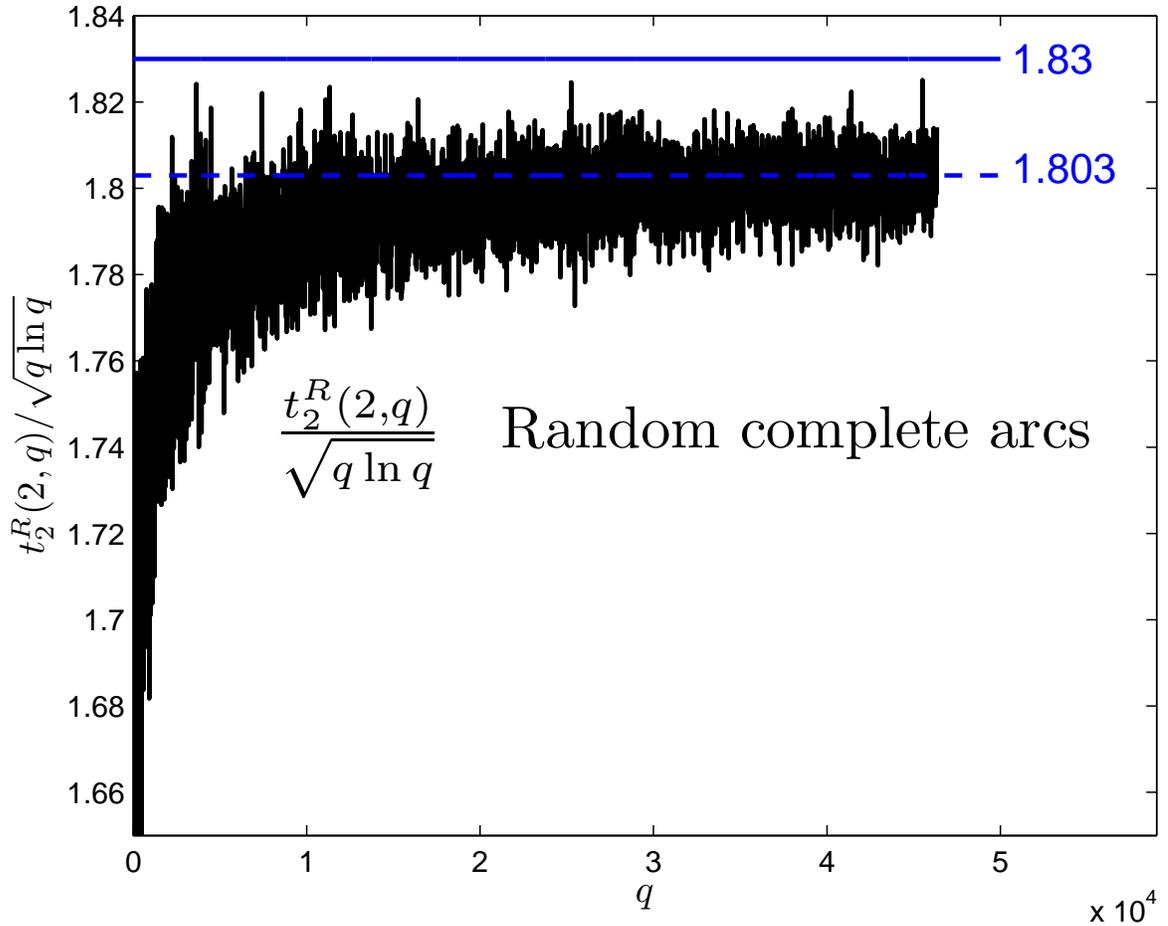,width=\textwidth}
\caption{\textbf{Bound $1.83\sqrt{q\ln q}$ vs random complete arcs.}
$y=1.83$ (the top solid line); $y=1.803$ (the 2-nd dashed line); values $t_{2}^{R}(2,q)/\sqrt{q\ln q}$,
$q\in R$,
where $t_{2}^{R}(2,q)$ is the size of a random complete arc (the solid curve)}
\label{fig1}
\end{figure}

In Figure \ref{fig1}, values $t_{2}^{R}(2,q)/\sqrt{q\ln q}$,
$q\in R$, are shown. The values oscillate around line
$y=1.803$; it means that the conjecture \eqref{eq1_aversize}
holds in the region $R$.

We denote
\begin{align}\label{eq2_Lex}
L^{\#}=\{195023,200003,205019,210011,215051,220009,225023,230003\}.
\end{align}
Sizes $t_{2}^{L}(2,q)$ of small complete arcs in $PG(2,q)$
obtained by Algorithm FOP with Lexicographical order in the
region $L^{\#}$ are obtained in this work. Values
$t_{2}^{L}(2,q)$ corresponding to $q$'s of \eqref{eq2_Lex} are
as follows:
\begin{align}\label{eq2_Mq}
\{2781,2822,2864,2886,2938,2958,3002,3033\}.
\end{align}

\begin{figure}[tbp]
\epsfig{file=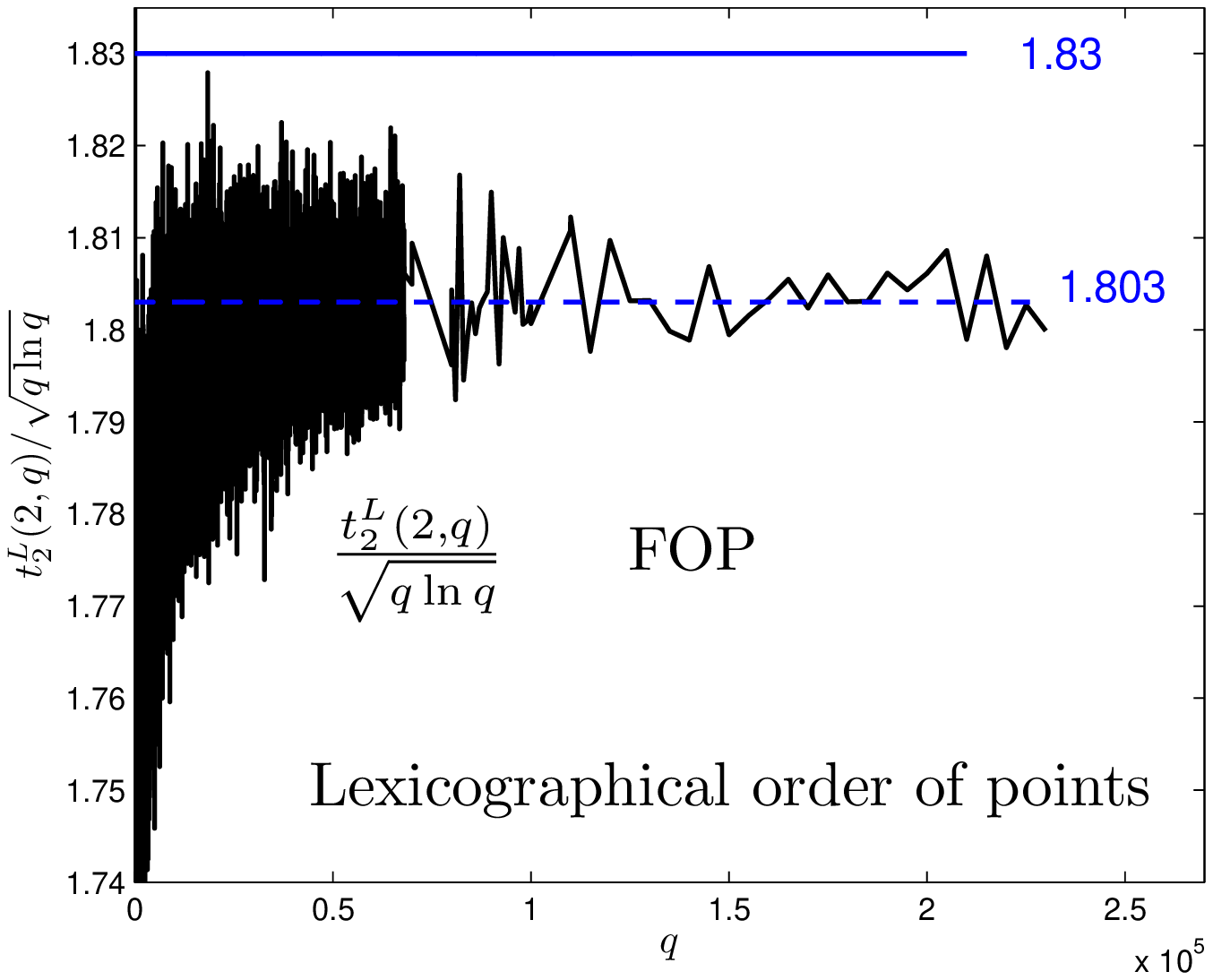,width=\textwidth}
\caption{\textbf{Bound $1.83\sqrt{q\ln q}$ vs Algorithm's FOP Lexicographical results.}
$y=1.83$ (the top solid line); $y=1.803$ (the 2-nd dashed line); values $t_{2}^{L}(2,q)/\sqrt{q\ln q}$, $q\in L\cup L^{\#}$,
where $t_{2}^{L}(2,q)$ is the size of a complete arc obtained by Algorithm FOP with Lexicographical order of points (the solid curve)}
\label{fig2}
\end{figure}

\begin{figure}[tbp]
\epsfig{file=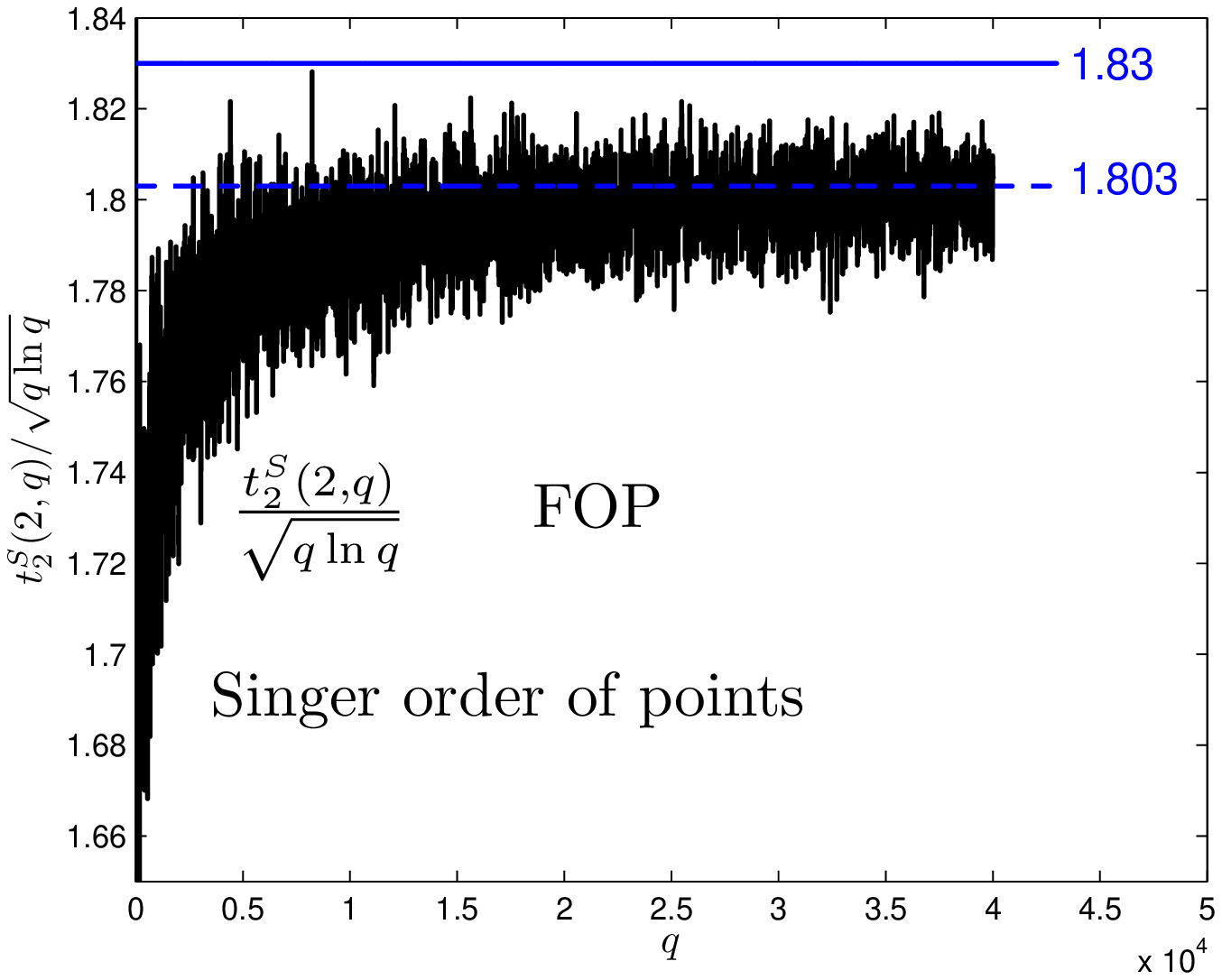,width=\textwidth}
\caption{\textbf{Bound $1.83\sqrt{q\ln q}$ vs Algorithm's FOP Singer results.}
$y=1.83$ (the top solid line); $y=1.803$ (the 2-nd dashed line); values $t_{2}^{S}(2,q)/\sqrt{q\ln q}$, $q\in S$,
where $t_{2}^{S}(2,q)$ is the size of a complete arc obtained by Algorithm FOP with Singer order of points (the solid curve)}
\label{fig3}
\end{figure}

In Figure \ref{fig2}, values $t_{2}^{L}(2,q)/\sqrt{q\ln q}$,
$q\in L\cup L^{\#}$, are shown. In Figure \ref{fig3}, values\\
$t_{2}^{S}(2,q)/\sqrt{q\ln q}$, $q\in S$, are given.

One can see that Figures \ref{fig1} and \ref{fig2}, \ref{fig3}
have the very similar structures. It is expected, as
Lexicographical order of points is a random order in the
geometrical sense. Singer order, of course, has a geometrical
sense but this sense is not connected with constructing of arcs
and with covering of points by bisecants. It is why Singer
order also may be treated as a random order.

\begin{figure}[tbp]
\epsfig{file=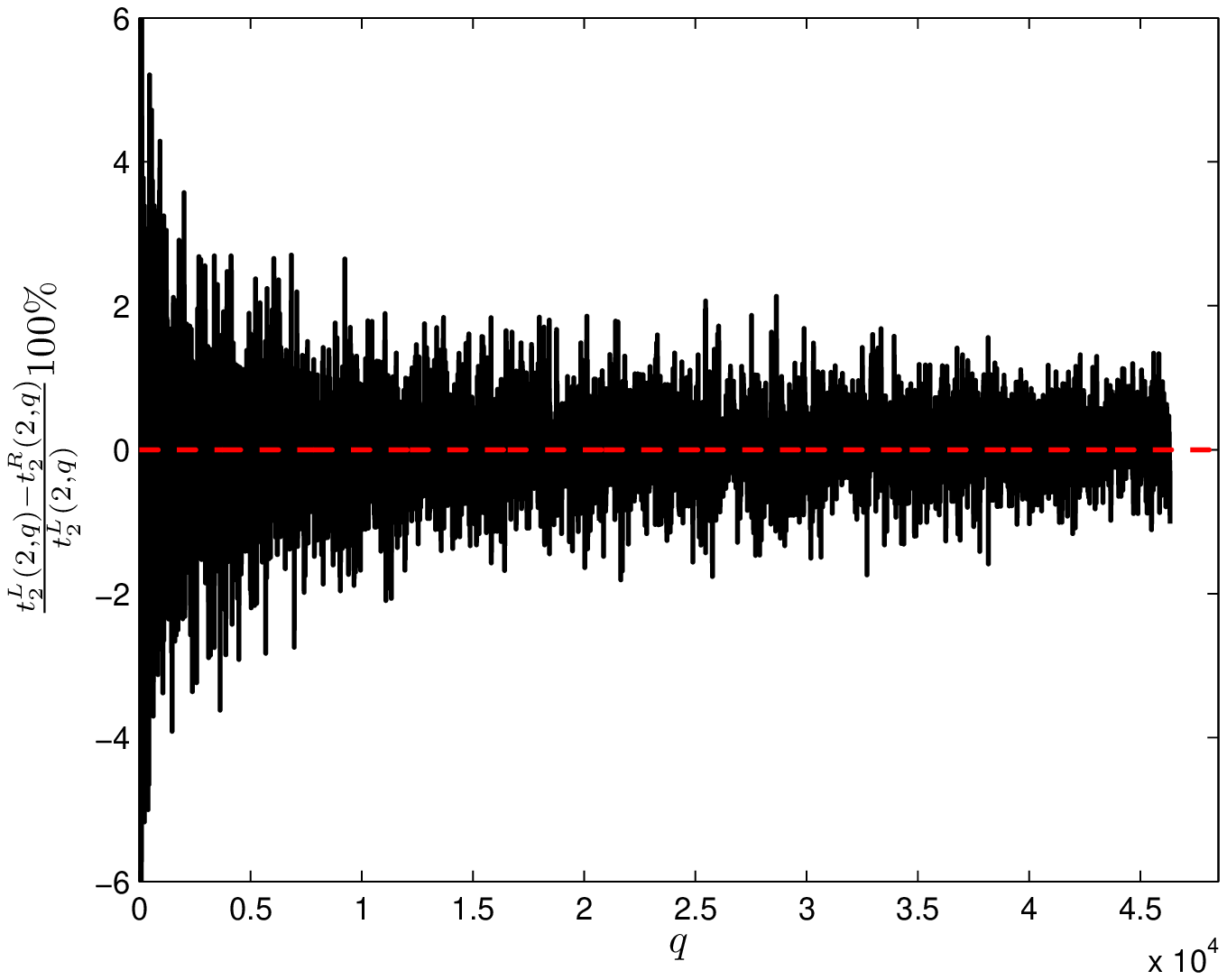,width=\textwidth}
\caption{\textbf{Difference $t^{L}_2(2,q)-t^{R}_2(2,q)$ in percentage.}
$\frac{t^{L}_2(2,q)-t^{R}_2(2,q)}{t^{L}_2(2,q)}100\%$, $q\in R$ (the solid curve)}
\label{fig4}
\end{figure}

\begin{figure}[tbp]
\epsfig{file=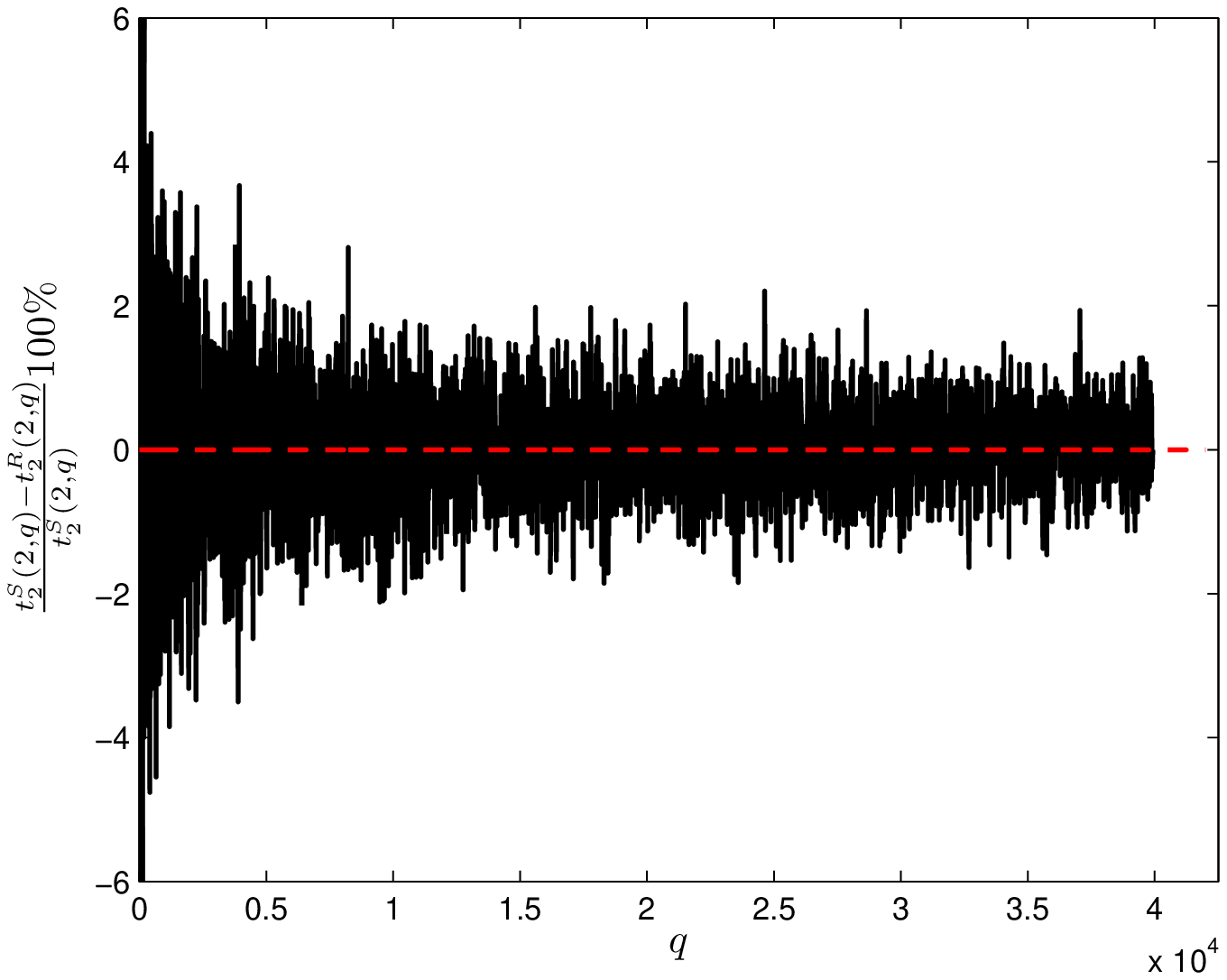,width=\textwidth}
\caption{\textbf{Difference $t^{S}_2(2,q)-t^{R}_2(2,q)$ in percentage.}
$\frac{t^{S}_2(2,q)-t^{R}_2(2,q)}{t^{S}_2(2,q)}100\%$, $q\in S$ (the solid curve)}
\label{fig5}
\end{figure}

In Figures \ref{fig4} and \ref{fig5}, the differences
$t^{L}_2(2,q)-t^{R}_2(2,q)$ and $t^{S}_2(2,q)-t^{S}_2(2,q)$ in
percentage are represented. We show values
$\frac{t^{L}_2(2,q)-t^{R}_2(2,q)}{t^{L}_2(2,q)}100\%$  in the
region $R$ and values of
$\frac{t^{S}_2(2,q)-t^{R}_2(2,q)}{t^{S}_2(2,q)}100\%$  in the
region $S$. Note that the differences  are less than 3\% for
$q>3000$.

\newpage
\textbf{Table 1} The sizes
 $t^{R}_{2}=t^{R}_{2}(2,q)$
 of random complete arcs in planes
 $\mathrm{PG}(2,q),$\\ $3\leq q\leq 46337$, $q$ prime  \medskip

 \renewcommand{\arraystretch}{1.00}
\begin{tabular}
 {@{}c@{\,\,\,}c|c@{\,\,\,}c|c@{\,\,\,}c|c@{\,\,\,}c|c@{\,\,\,}c|c@{\,\,\,}c@{}}
  \hline
$q$ & $t^{R}_{2}$ &
$q$ & $t^{R}_{2}$ &
$q$ & $t^{R}_{2}$ &
$q$ & $t^{R}_{2}$ &
$q$ & $t^{R}_{2}$ &
$q$ & $t^{R}_{2}\vphantom{H^{H^{L^{L}}}}$ \\
 \hline
 3 & 4  &  5 & 6  &  7 & 6  &  11 & 8  &  13 & 10  &  17 & 11  \\
 19 & 13  &  23 & 13  &  29 & 15  &  31 & 16  &  37 & 19  &  41 & 20  \\
 43 & 20  &  47 & 21  &  53 & 24  &  59 & 24  &  61 & 25  &  67 & 28  \\
 71 & 30  &  73 & 28  &  79 & 30  &  83 & 33  &  89 & 34  &  97 & 33  \\
 101 & 37  &  103 & 35  &  107 & 35  &  109 & 36  &  113 & 38  &  127 & 42  \\
 131 & 41  &  137 & 45  &  139 & 44  &  149 & 46  &  151 & 46  &  157 & 47  \\
 163 & 47  &  167 & 48  &  173 & 49  &  179 & 51  &  181 & 52  &  191 & 52  \\
\hline
 \end{tabular}

 \newpage

\textbf{Table 1} Continue 1 \medskip

 \renewcommand{\arraystretch}{1.00}
 \begin{tabular}
 {@{}c@{\,\,\,}c|c@{\,\,\,}c|c@{\,\,\,}c|c@{\,\,\,}c|c@{\,\,\,}c|c@{\,\,\,}c@{}}
  \hline
$q$ & $t^{R}_{2}$ &
$q$ & $t^{R}_{2}$ &
$q$ & $t^{R}_{2}$ &
$q$ & $t^{R}_{2}$ &
$q$ & $t^{R}_{2}$ &
$q$ & $t^{R}_{2}\vphantom{H^{H^{L^{L}}}}$ \\
 \hline
 1301 & 170  &  1303 & 170  &  1307 & 170  &  1319 & 174  &  1321 & 170  &  1327 & 173  \\
 1361 & 175  &  1367 & 174  &  1373 & 173  &  1381 & 174  &  1399 & 180  &  1409 & 177  \\
 1423 & 178  &  1427 & 178  &  1429 & 176  &  1433 & 180  &  1439 & 177  &  1447 & 182  \\
 1451 & 183  &  1453 & 179  &  1459 & 181  &  1471 & 186  &  1481 & 184  &  1483 & 185  \\
 1487 & 182  &  1489 & 181  &  1493 & 185  &  1499 & 185  &  1511 & 183  &  1523 & 183  \\
 1531 & 185  &  1543 & 187  &  1549 & 187  &  1553 & 189  &  1559 & 187  &  1567 & 187  \\
 1571 & 188  &  1579 & 187  &  1583 & 189  &  1597 & 193  &  1601 & 192  &  1607 & 193  \\
\hline
 \end{tabular}

 \newpage

\textbf{Table 1} Continue 2 \medskip

 \renewcommand{\arraystretch}{1.00}
 \begin{tabular}
 {@{}c@{\,\,\,}c|c@{\,\,\,}c|c@{\,\,\,}c|c@{\,\,\,}c|c@{\,\,\,}c|c@{\,\,\,}c@{}}
  \hline
$q$ & $t^{R}_{2}$ &
$q$ & $t^{R}_{2}$ &
$q$ & $t^{R}_{2}$ &
$q$ & $t^{R}_{2}$ &
$q$ & $t^{R}_{2}$ &
$q$ & $t^{R}_{2}\vphantom{H^{H^{L^{L}}}}$ \\
 \hline
 2917 & 268  &  2927 & 271  &  2939 & 271  &  2953 & 272  &  2957 & 267  &  2963 & 269  \\
 2969 & 275  &  2971 & 274  &  2999 & 274  &  3001 & 272  &  3011 & 275  &  3019 & 274  \\
 3023 & 279  &  3037 & 276  &  3041 & 272  &  3049 & 281  &  3061 & 277  &  3067 & 276  \\
 3079 & 275  &  3083 & 279  &  3089 & 281  &  3109 & 279  &  3119 & 285  &  3121 & 283  \\
 3137 & 285  &  3163 & 279  &  3167 & 283  &  3169 & 282  &  3181 & 283  &  3187 & 285  \\
 3191 & 280  &  3203 & 288  &  3209 & 286  &  3217 & 283  &  3221 & 288  &  3229 & 288  \\
 3251 & 287  &  3253 & 290  &  3257 & 289  &  3259 & 289  &  3271 & 290  &  3299 & 290  \\
\hline
 \end{tabular}

 \newpage

\textbf{Table 1} Continue 3 \medskip

 \renewcommand{\arraystretch}{1.00}
 \begin{tabular}
 {@{}c@{\,\,\,}c|c@{\,\,\,}c|c@{\,\,\,}c|c@{\,\,\,}c|c@{\,\,\,}c|c@{\,\,\,}c@{}}
  \hline
$q$ & $t^{R}_{2}$ &
$q$ & $t^{R}_{2}$ &
$q$ & $t^{R}_{2}$ &
$q$ & $t^{R}_{2}$ &
$q$ & $\overline{t}_{2}y$ &
$q$ & $t^{R}_{2}\vphantom{H^{H^{L^{L}}}}$ \\
 \hline
 4673 & 354  &  4679 & 355  &  4691 & 357  &  4703 & 357  &  4721 & 358  &  4723 & 354  \\
 4729 & 354  &  4733 & 356  &  4751 & 357  &  4759 & 360  &  4783 & 361  &  4787 & 357  \\
 4789 & 362  &  4793 & 357  &  4799 & 362  &  4801 & 360  &  4813 & 360  &  4817 & 361  \\
 4831 & 357  &  4861 & 362  &  4871 & 364  &  4877 & 361  &  4889 & 364  &  4903 & 363  \\
 4909 & 364  &  4919 & 365  &  4931 & 362  &  4933 & 362  &  4937 & 365  &  4943 & 363  \\
 4951 & 367  &  4957 & 363  &  4967 & 368  &  4969 & 368  &  4973 & 365  &  4987 & 365  \\
 4993 & 367  &  4999 & 366  &  5003 & 370  &  5009 & 366  &  5011 & 370  &  5021 & 371  \\
\hline
 \end{tabular}

 \newpage

\textbf{Table 1} Continue 4 \medskip

 \renewcommand{\arraystretch}{1.00}
 \begin{tabular}
 {@{}c@{\,\,\,}c|c@{\,\,\,}c|c@{\,\,\,}c|c@{\,\,\,}c|c@{\,\,\,}c|c@{\,\,\,}c@{}}
  \hline
$q$ & $t^{R}_{2}$ &
$q$ & $t^{R}_{2}$ &
$q$ & $t^{R}_{2}$ &
$q$ & $t^{R}_{2}$ &
$q$ & $t^{R}_{2}$ &
$q$ & $t^{R}_{2}\vphantom{H^{H^{L^{L}}}}$ \\
 \hline
 6491 & 427  &  6521 & 423  &  6529 & 426  &  6547 & 427  &  6551 & 425  &  6553 & 432  \\
 6563 & 423  &  6569 & 431  &  6571 & 429  &  6577 & 432  &  6581 & 425  &  6599 & 427  \\
 6607 & 434  &  6619 & 432  &  6637 & 431  &  6653 & 429  &  6659 & 435  &  6661 & 429  \\
 6673 & 435  &  6679 & 431  &  6689 & 435  &  6691 & 434  &  6701 & 434  &  6703 & 428  \\
 6709 & 434  &  6719 & 435  &  6733 & 438  &  6737 & 437  &  6761 & 435  &  6763 & 441  \\
 6779 & 439  &  6781 & 436  &  6791 & 435  &  6793 & 437  &  6803 & 433  &  6823 & 436  \\
 6827 & 437  &  6829 & 436  &  6833 & 432  &  6841 & 437  &  6857 & 441  &  6863 & 440  \\
\hline
 \end{tabular}

 \newpage

\textbf{Table 1} Continue 5 \medskip

 \renewcommand{\arraystretch}{1.00}
 \begin{tabular}
 {@{}c@{\,\,\,}c|c@{\,\,\,}c|c@{\,\,\,}c|c@{\,\,\,}c|c@{\,\,\,}c|c@{\,\,\,}c@{}}
  \hline
$q$ & $t^{R}_{2}$ &
$q$ & $t^{R}_{2}$ &
$q$ & $t^{R}_{2}$ &
$q$ & $t^{R}_{2}$ &
$q$ & $t^{R}_{2}$ &
$q$ & $t^{R}_{2}\vphantom{H^{H^{L^{L}}}}$ \\
 \hline
 8419 & 496  &  8423 & 495  &  8429 & 496  &  8431 & 495  &  8443 & 494  &  8447 & 494  \\
 8461 & 497  &  8467 & 497  &  8501 & 495  &  8513 & 497  &  8521 & 495  &  8527 & 498  \\
 8537 & 493  &  8539 & 497  &  8543 & 498  &  8563 & 495  &  8573 & 497  &  8581 & 495  \\
 8597 & 500  &  8599 & 497  &  8609 & 502  &  8623 & 496  &  8627 & 499  &  8629 & 497  \\
 8641 & 505  &  8647 & 504  &  8663 & 507  &  8669 & 503  &  8677 & 500  &  8681 & 501  \\
 8689 & 503  &  8693 & 507  &  8699 & 506  &  8707 & 504  &  8713 & 507  &  8719 & 505  \\
 8731 & 506  &  8737 & 502  &  8741 & 506  &  8747 & 504  &  8753 & 508  &  8761 & 501  \\
\hline
 \end{tabular}

 \newpage

\textbf{Table 1} Continue 6 \medskip

 \renewcommand{\arraystretch}{1.00}
 \begin{tabular}
 {@{}c@{\,\,\,}c|c@{\,\,\,}c|c@{\,\,\,}c|c@{\,\,\,}c|c@{\,\,\,}c|c@{\,\,\,}c@{}}
  \hline
$q$ & $t^{R}_{2}$ &
$q$ & $t^{R}_{2}$ &
$q$ & $t^{R}_{2}$ &
$q$ & $t^{R}_{2}$ &
$q$ & $t^{R}_{2}$ &
$q$ & $t^{R}_{2}\vphantom{H^{H^{L^{L}}}}$ \\
 \hline
 10273 & 549  &  10289 & 549  &  10301 & 555  &  10303 & 557  &  10313 & 550  &  10321 & 553  \\
 10331 & 555  &  10333 & 554  &  10337 & 554  &  10343 & 555  &  10357 & 556  &  10369 & 556  \\
 10391 & 554  &  10399 & 553  &  10427 & 555  &  10429 & 554  &  10433 & 552  &  10453 & 564  \\
 10457 & 555  &  10459 & 551  &  10463 & 556  &  10477 & 558  &  10487 & 560  &  10499 & 557  \\
 10501 & 559  &  10513 & 558  &  10529 & 559  &  10531 & 562  &  10559 & 561  &  10567 & 564  \\
 10589 & 563  &  10597 & 563  &  10601 & 564  &  10607 & 564  &  10613 & 563  &  10627 & 564  \\
 10631 & 557  &  10639 & 563  &  10651 & 560  &  10657 & 565  &  10663 & 567  &  10667 & 561  \\
\hline
 \end{tabular}

 \newpage

\textbf{Table 1} Continue 7 \medskip

 \renewcommand{\arraystretch}{1.00}
 \begin{tabular}
 {@{}c@{\,\,\,}c|c@{\,\,\,}c|c@{\,\,\,}c|c@{\,\,\,}c|c@{\,\,\,}c|c@{\,\,\,}c@{}}
  \hline
$q$ & $t^{R}_{2}$ &
$q$ & $t^{R}_{2}$ &
$q$ & $t^{R}_{2}$ &
$q$ & $t^{R}_{2}$ &
$q$ & $\overline{t}_{2}y$ &
$q$ & $t^{R}_{2}\vphantom{H^{H^{L^{L}}}}$ \\
 \hline
 12301 & 613  &  12323 & 608  &  12329 & 607  &  12343 & 616  &  12347 & 609  &  12373 & 619  \\
 12377 & 613  &  12379 & 611  &  12391 & 617  &  12401 & 612  &  12409 & 610  &  12413 & 613  \\
 12421 & 612  &  12433 & 610  &  12437 & 613  &  12451 & 616  &  12457 & 608  &  12473 & 618  \\
 12479 & 612  &  12487 & 617  &  12491 & 614  &  12497 & 620  &  12503 & 617  &  12511 & 618  \\
 12517 & 614  &  12527 & 616  &  12539 & 622  &  12541 & 620  &  12547 & 617  &  12553 & 617  \\
 12569 & 624  &  12577 & 612  &  12583 & 618  &  12589 & 618  &  12601 & 622  &  12611 & 619  \\
 12613 & 621  &  12619 & 619  &  12637 & 622  &  12641 & 622  &  12647 & 628  &  12653 & 618  \\
\hline
 \end{tabular}

 \newpage

\textbf{Table 1} Continue 8 \medskip

 \renewcommand{\arraystretch}{1.00}
 \begin{tabular}
 {@{}c@{\,\,\,}c|c@{\,\,\,}c|c@{\,\,\,}c|c@{\,\,\,}c|c@{\,\,\,}c|c@{\,\,\,}c@{}}
  \hline
$q$ & $t^{R}_{2}$ &
$q$ & $t^{R}_{2}$ &
$q$ & $t^{R}_{2}$ &
$q$ & $t^{R}_{2}$ &
$q$ & $t^{R}_{2}$ &
$q$ & $t^{R}_{2}\vphantom{H^{H^{L^{L}}}}$ \\
 \hline
 14327 & 660  &  14341 & 666  &  14347 & 667  &  14369 & 663  &  14387 & 664  &  14389 & 661  \\
 14401 & 664  &  14407 & 669  &  14411 & 665  &  14419 & 668  &  14423 & 670  &  14431 & 665  \\
 14437 & 668  &  14447 & 668  &  14449 & 670  &  14461 & 671  &  14479 & 667  &  14489 & 672  \\
 14503 & 667  &  14519 & 676  &  14533 & 671  &  14537 & 668  &  14543 & 666  &  14549 & 668  \\
 14551 & 665  &  14557 & 667  &  14561 & 668  &  14563 & 668  &  14591 & 674  &  14593 & 670  \\
 14621 & 672  &  14627 & 671  &  14629 & 675  &  14633 & 671  &  14639 & 672  &  14653 & 674  \\
 14657 & 676  &  14669 & 676  &  14683 & 666  &  14699 & 671  &  14713 & 672  &  14717 & 677  \\
\hline
 \end{tabular}

 \newpage

\textbf{Table 1} Continue 9 \medskip

 \renewcommand{\arraystretch}{1.00}
 \begin{tabular}
 {@{}c@{\,\,\,}c|c@{\,\,\,}c|c@{\,\,\,}c|c@{\,\,\,}c|c@{\,\,\,}c|c@{\,\,\,}c@{}}
  \hline
$q$ & $t^{R}_{2}$ &
$q$ & $t^{R}_{2}$ &
$q$ & $t^{R}_{2}$ &
$q$ & $t^{R}_{2}$ &
$q$ & $t^{R}_{2}$ &
$q$ & $t^{R}_{2}\vphantom{H^{H^{L^{L}}}}$ \\
 \hline
 16267 & 715  &  16273 & 708  &  16301 & 710  &  16319 & 719  &  16333 & 709  &  16339 & 713  \\
 16349 & 714  &  16361 & 718  &  16363 & 719  &  16369 & 715  &  16381 & 712  &  16411 & 717  \\
 16417 & 719  &  16421 & 717  &  16427 & 727  &  16433 & 717  &  16447 & 716  &  16451 & 716  \\
 16453 & 717  &  16477 & 714  &  16481 & 722  &  16487 & 721  &  16493 & 714  &  16519 & 719  \\
 16529 & 719  &  16547 & 719  &  16553 & 722  &  16561 & 723  &  16567 & 717  &  16573 & 722  \\
 16603 & 717  &  16607 & 720  &  16619 & 718  &  16631 & 716  &  16633 & 719  &  16649 & 718  \\
 16651 & 722  &  16657 & 720  &  16661 & 724  &  16673 & 725  &  16691 & 719  &  16693 & 725  \\
\hline
 \end{tabular}

 \newpage

\textbf{Table 1} Continue 10 \medskip

 \renewcommand{\arraystretch}{1.00}
 \begin{tabular}
 {@{}c@{\,\,\,}c|c@{\,\,\,}c|c@{\,\,\,}c|c@{\,\,\,}c|c@{\,\,\,}c|c@{\,\,\,}c@{}}
  \hline
$q$ & $t^{R}_{2}$ &
$q$ & $t^{R}_{2}$ &
$q$ & $t^{R}_{2}$ &
$q$ & $t^{R}_{2}$ &
$q$ & $t^{R}_{2}$ &
$q$ & $t^{R}_{2}\vphantom{H^{H^{L^{L}}}}$ \\
 \hline
 18313 & 761  &  18329 & 762  &  18341 & 758  &  18353 & 768  &  18367 & 765  &  18371 & 765  \\
 18379 & 763  &  18397 & 762  &  18401 & 760  &  18413 & 759  &  18427 & 764  &  18433 & 764  \\
 18439 & 772  &  18443 & 764  &  18451 & 768  &  18457 & 765  &  18461 & 772  &  18481 & 762  \\
 18493 & 770  &  18503 & 767  &  18517 & 767  &  18521 & 767  &  18523 & 768  &  18539 & 768  \\
 18541 & 770  &  18553 & 767  &  18583 & 768  &  18587 & 767  &  18593 & 769  &  18617 & 764  \\
 18637 & 774  &  18661 & 772  &  18671 & 766  &  18679 & 772  &  18691 & 778  &  18701 & 774  \\
 18713 & 774  &  18719 & 770  &  18731 & 769  &  18743 & 772  &  18749 & 766  &  18757 & 765  \\
\hline
 \end{tabular}

 \newpage

\textbf{Table 1} Continue 11 \medskip

 \renewcommand{\arraystretch}{1.00}
 \begin{tabular}
 {@{}c@{\,\,\,}c|c@{\,\,\,}c|c@{\,\,\,}c|c@{\,\,\,}c|c@{\,\,\,}c|c@{\,\,\,}c@{}}
  \hline
$q$ & $t^{R}_{2}$ &
$q$ & $t^{R}_{2}$ &
$q$ & $t^{R}_{2}$ &
$q$ & $t^{R}_{2}$ &
$q$ & $t^{R}_{2}$ &
$q$ & $t^{R}_{2}\vphantom{H^{H^{L^{L}}}}$ \\
 \hline
 20441 & 808  &  20443 & 806  &  20477 & 808  &  20479 & 815  &  20483 & 807  &  20507 & 807  \\
 20509 & 813  &  20521 & 810  &  20533 & 815  &  20543 & 815  &  20549 & 812  &  20551 & 814  \\
 20563 & 813  &  20593 & 813  &  20599 & 815  &  20611 & 816  &  20627 & 814  &  20639 & 812  \\
 20641 & 819  &  20663 & 813  &  20681 & 819  &  20693 & 811  &  20707 & 812  &  20717 & 814  \\
 20719 & 812  &  20731 & 818  &  20743 & 819  &  20747 & 817  &  20749 & 811  &  20753 & 815  \\
 20759 & 815  &  20771 & 818  &  20773 & 819  &  20789 & 816  &  20807 & 817  &  20809 & 812  \\
 20849 & 817  &  20857 & 821  &  20873 & 822  &  20879 & 820  &  20887 & 819  &  20897 & 824  \\
\hline
 \end{tabular}

 \newpage

\textbf{Table 1} Continue 12 \medskip

 \renewcommand{\arraystretch}{1.00}
 \begin{tabular}
 {@{}c@{\,\,\,}c|c@{\,\,\,}c|c@{\,\,\,}c|c@{\,\,\,}c|c@{\,\,\,}c|c@{\,\,\,}c@{}}
  \hline
$q$ & $t^{R}_{2}$ &
$q$ & $t^{R}_{2}$ &
$q$ & $t^{R}_{2}$ &
$q$ & $t^{R}_{2}$ &
$q$ & $t^{R}_{2}$ &
$q$ & $t^{R}_{2}\vphantom{H^{H^{L^{L}}}}$ \\
 \hline
 22541 & 853  &  22543 & 856  &  22549 & 857  &  22567 & 858  &  22571 & 859  &  22573 & 856  \\
 22613 & 853  &  22619 & 854  &  22621 & 857  &  22637 & 858  &  22639 & 857  &  22643 & 853  \\
 22651 & 855  &  22669 & 863  &  22679 & 857  &  22691 & 864  &  22697 & 859  &  22699 & 860  \\
 22709 & 857  &  22717 & 856  &  22721 & 852  &  22727 & 856  &  22739 & 858  &  22741 & 864  \\
 22751 & 856  &  22769 & 853  &  22777 & 863  &  22783 & 853  &  22787 & 856  &  22807 & 857  \\
 22811 & 868  &  22817 & 861  &  22853 & 858  &  22859 & 856  &  22861 & 859  &  22871 & 866  \\
 22877 & 866  &  22901 & 857  &  22907 & 863  &  22921 & 867  &  22937 & 859  &  22943 & 866  \\
\hline
 \end{tabular}

 \newpage

\textbf{Table 1} Continue 13 \medskip

 \renewcommand{\arraystretch}{1.00}
 \begin{tabular}
 {@{}c@{\,\,\,}c|c@{\,\,\,}c|c@{\,\,\,}c|c@{\,\,\,}c|c@{\,\,\,}c|c@{\,\,\,}c@{}}
  \hline
$q$ & $t^{R}_{2}$ &
$q$ & $t^{R}_{2}$ &
$q$ & $t^{R}_{2}$ &
$q$ & $t^{R}_{2}$ &
$q$ & $t^{R}_{2}$ &
$q$ & $t^{R}_{2}\vphantom{H^{H^{L^{L}}}}$ \\
 \hline
 24631 & 888  &  24659 & 901  &  24671 & 896  &  24677 & 899  &  24683 & 898  &  24691 & 904  \\
 24697 & 899  &  24709 & 903  &  24733 & 899  &  24749 & 895  &  24763 & 894  &  24767 & 897  \\
 24781 & 900  &  24793 & 906  &  24799 & 896  &  24809 & 901  &  24821 & 904  &  24841 & 903  \\
 24847 & 908  &  24851 & 901  &  24859 & 907  &  24877 & 906  &  24889 & 910  &  24907 & 902  \\
 24917 & 902  &  24919 & 905  &  24923 & 902  &  24943 & 905  &  24953 & 903  &  24967 & 903  \\
 24971 & 903  &  24977 & 905  &  24979 & 905  &  24989 & 901  &  24989 & 911  &  25013 & 907  \\
 25031 & 909  &  25033 & 900  &  25037 & 908  &  25057 & 901  &  25073 & 905  &  25087 & 910  \\
\hline
 \end{tabular}

 \newpage

\textbf{Table 1} Continue 14 \medskip

 \renewcommand{\arraystretch}{1.00}
 \begin{tabular}
 {@{}c@{\,\,\,}c|c@{\,\,\,}c|c@{\,\,\,}c|c@{\,\,\,}c|c@{\,\,\,}c|c@{\,\,\,}c@{}}
  \hline
$q$ & $t^{R}_{2}$ &
$q$ & $t^{R}_{2}$ &
$q$ & $t^{R}_{2}$ &
$q$ & $t^{R}_{2}$ &
$q$ & $t^{R}_{2}$ &
$q$ & $t^{R}_{2}\vphantom{H^{H^{L^{L}}}}$ \\
 \hline
 26777 & 935  &  26783 & 935  &  26801 & 939  &  26813 & 939  &  26821 & 945  &  26833 & 937  \\
 26839 & 944  &  26849 & 936  &  26861 & 937  &  26863 & 943  &  26879 & 940  &  26881 & 942  \\
 26891 & 947  &  26893 & 946  &  26903 & 937  &  26921 & 941  &  26927 & 940  &  26947 & 947  \\
 26951 & 942  &  26953 & 939  &  26959 & 941  &  26981 & 940  &  26987 & 942  &  26993 & 940  \\
 27011 & 954  &  27017 & 942  &  27031 & 945  &  27043 & 951  &  27059 & 944  &  27061 & 942  \\
 27067 & 944  &  27073 & 950  &  27077 & 939  &  27091 & 944  &  27103 & 945  &  27107 & 949  \\
 27109 & 945  &  27127 & 949  &  27143 & 944  &  27179 & 946  &  27191 & 945  &  27197 & 950  \\
\hline
 \end{tabular}

 \newpage

\textbf{Table 1} Continue 15 \medskip

 \renewcommand{\arraystretch}{1.00}
 \begin{tabular}
 {@{}c@{\,\,\,}c|c@{\,\,\,}c|c@{\,\,\,}c|c@{\,\,\,}c|c@{\,\,\,}c|c@{\,\,\,}c@{}}
  \hline
$q$ & $t^{R}_{2}$ &
$q$ & $t^{R}_{2}$ &
$q$ & $t^{R}_{2}$ &
$q$ & $t^{R}_{2}$ &
$q$ & $t^{R}_{2}$ &
$q$ & $t^{R}_{2}\vphantom{H^{H^{L^{L}}}}$ \\
 \hline
 28921 & 980  &  28927 & 978  &  28933 & 982  &  28949 & 976  &  28961 & 985  &  28979 & 982  \\
 29009 & 988  &  29017 & 985  &  29021 & 981  &  29023 & 980  &  29027 & 978  &  29033 & 980  \\
 29059 & 983  &  29063 & 980  &  29077 & 989  &  29101 & 985  &  29123 & 986  &  29129 & 985  \\
 29131 & 984  &  29137 & 993  &  29147 & 988  &  29153 & 987  &  29167 & 987  &  29173 & 989  \\
 29179 & 988  &  29191 & 982  &  29201 & 996  &  29207 & 983  &  29209 & 992  &  29221 & 986  \\
 29231 & 985  &  29243 & 985  &  29251 & 991  &  29269 & 985  &  29287 & 986  &  29297 & 989  \\
 29303 & 998  &  29311 & 994  &  29327 & 990  &  29333 & 985  &  29339 & 986  &  29347 & 989  \\
\hline
 \end{tabular}

 \newpage

\textbf{Table 1} Continue 16 \medskip

 \renewcommand{\arraystretch}{1.00}
 \begin{tabular}
 {@{}c@{\,\,\,}c|c@{\,\,\,}c|c@{\,\,\,}c|c@{\,\,\,}c|c@{\,\,\,}c|c@{\,\,\,}c@{}}
  \hline
$q$ & $t^{R}_{2}$ &
$q$ & $t^{R}_{2}$ &
$q$ & $t^{R}_{2}$ &
$q$ & $t^{R}_{2}$ &
$q$ & $t^{R}_{2}$ &
$q$ & $t^{R}_{2}\vphantom{H^{H^{L^{L}}}}$ \\
 \hline
 31159 & 1021  &  31177 & 1026  &  31181 & 1020  &  31183 & 1022  &  31189 & 1017  &  31193 & 1026  \\
 31219 & 1025  &  31223 & 1021  &  31231 & 1018  &  31237 & 1027  &  31247 & 1024  &  31249 & 1019  \\
 31253 & 1023  &  31259 & 1025  &  31267 & 1026  &  31271 & 1022  &  31277 & 1018  &  31307 & 1028  \\
 31319 & 1020  &  31321 & 1027  &  31327 & 1030  &  31333 & 1030  &  31337 & 1024  &  31357 & 1028  \\
 31379 & 1032  &  31387 & 1019  &  31391 & 1031  &  31393 & 1031  &  31397 & 1030  &  31469 & 1025  \\
 31477 & 1023  &  31481 & 1031  &  31489 & 1028  &  31511 & 1028  &  31513 & 1029  &  31517 & 1033  \\
 31531 & 1024  &  31541 & 1031  &  31543 & 1028  &  31547 & 1025  &  31567 & 1024  &  31573 & 1030  \\
\hline
 \end{tabular}

 \newpage

\textbf{Table 1} Continue 17 \medskip

 \renewcommand{\arraystretch}{1.00}
 \begin{tabular}
 {@{}c@{\,\,\,}c|c@{\,\,\,}c|c@{\,\,\,}c|c@{\,\,\,}c|c@{\,\,\,}c|c@{\,\,\,}c@{}}
  \hline
$q$ & $t^{R}_{2}$ &
$q$ & $t^{R}_{2}$ &
$q$ & $t^{R}_{2}$ &
$q$ & $t^{R}_{2}$ &
$q$ & $t^{R}_{2}$ &
$q$ & $t^{R}_{2}\vphantom{H^{H^{L^{L}}}}$ \\
 \hline
 33311 & 1064  &  33317 & 1057  &  33329 & 1064  &  33331 & 1058  &  33343 & 1056  &  33347 & 1064  \\
 33349 & 1053  &  33353 & 1061  &  33359 & 1065  &  33377 & 1067  &  33391 & 1063  &  33403 & 1059  \\
 33409 & 1058  &  33413 & 1062  &  33427 & 1065  &  33457 & 1061  &  33461 & 1062  &  33469 & 1058  \\
 33479 & 1063  &  33487 & 1060  &  33493 & 1074  &  33503 & 1064  &  33521 & 1059  &  33529 & 1068  \\
 33533 & 1063  &  33547 & 1058  &  33563 & 1067  &  33569 & 1069  &  33577 & 1071  &  33581 & 1065  \\
 33587 & 1073  &  33589 & 1064  &  33599 & 1063  &  33601 & 1070  &  33613 & 1071  &  33617 & 1069  \\
 33619 & 1061  &  33623 & 1069  &  33629 & 1069  &  33637 & 1069  &  33641 & 1068  &  33647 & 1067  \\
\hline
 \end{tabular}

 \newpage

\textbf{Table 1} Continue 18 \medskip

 \renewcommand{\arraystretch}{1.00}
 \begin{tabular}
 {@{}c@{\,\,\,}c|c@{\,\,\,}c|c@{\,\,\,}c|c@{\,\,\,}c|c@{\,\,\,}c|c@{\,\,\,}c@{}}
  \hline
$q$ & $t^{R}_{2}$ &
$q$ & $t^{R}_{2}$ &
$q$ & $t^{R}_{2}$ &
$q$ & $t^{R}_{2}$ &
$q$ & $t^{R}_{2}$ &
$q$ & $t^{R}_{2}\vphantom{H^{H^{L^{L}}}}$ \\
 \hline
 35449 & 1096  &  35461 & 1095  &  35491 & 1094  &  35507 & 1098  &  35509 & 1097  &  35521 & 1096  \\
 35527 & 1098  &  35531 & 1095  &  35533 & 1099  &  35537 & 1093  &  35543 & 1097  &  35569 & 1104  \\
 35573 & 1102  &  35591 & 1108  &  35593 & 1104  &  35597 & 1102  &  35603 & 1100  &  35617 & 1099  \\
 35671 & 1105  &  35677 & 1096  &  35729 & 1108  &  35731 & 1099  &  35747 & 1107  &  35753 & 1111  \\
 35759 & 1102  &  35771 & 1099  &  35797 & 1111  &  35801 & 1099  &  35803 & 1111  &  35809 & 1104  \\
 35831 & 1106  &  35837 & 1107  &  35839 & 1100  &  35851 & 1105  &  35863 & 1103  &  35869 & 1100  \\
 35879 & 1100  &  35897 & 1104  &  35899 & 1099  &  35911 & 1098  &  35923 & 1106  &  35933 & 1113  \\
\hline
 \end{tabular}

 \newpage

\textbf{Table 1} Continue 19 \medskip

 \renewcommand{\arraystretch}{1.00}
 \begin{tabular}
 {@{}c@{\,\,\,}c|c@{\,\,\,}c|c@{\,\,\,}c|c@{\,\,\,}c|c@{\,\,\,}c|c@{\,\,\,}c@{}}
  \hline
$q$ & $t^{R}_{2}$ &
$q$ & $t^{R}_{2}$ &
$q$ & $t^{R}_{2}$ &
$q$ & $t^{R}_{2}$ &
$q$ & $t^{R}_{2}$ &
$q$ & $t^{R}_{2}\vphantom{H^{H^{L^{L}}}}$ \\
 \hline
 37619 & 1128  &  37633 & 1132  &  37643 & 1137  &  37649 & 1133  &  37657 & 1133  &  37663 & 1129  \\
 37691 & 1130  &  37693 & 1131  &  37699 & 1132  &  37717 & 1135  &  37747 & 1142  &  37781 & 1135  \\
 37783 & 1147  &  37799 & 1135  &  37811 & 1139  &  37813 & 1138  &  37831 & 1142  &  37847 & 1142  \\
 37853 & 1138  &  37861 & 1139  &  37871 & 1139  &  37879 & 1144  &  37889 & 1134  &  37897 & 1139  \\
 37907 & 1138  &  37951 & 1131  &  37957 & 1135  &  37963 & 1137  &  37967 & 1143  &  37987 & 1128  \\
 37991 & 1136  &  37993 & 1151  &  37997 & 1137  &  38011 & 1136  &  38039 & 1139  &  38047 & 1142  \\
 38053 & 1141  &  38069 & 1146  &  38083 & 1144  &  38113 & 1141  &  38119 & 1139  &  38149 & 1139  \\
\hline
 \end{tabular}

 \newpage

\textbf{Table 1} Continue 20 \medskip

 \renewcommand{\arraystretch}{1.00}
 \begin{tabular}
 {@{}c@{\,\,\,}c|c@{\,\,\,}c|c@{\,\,\,}c|c@{\,\,\,}c|c@{\,\,\,}c|c@{\,\,\,}c@{}}
  \hline
$q$ & $t^{R}_{2}$ &
$q$ & $t^{R}_{2}$ &
$q$ & $t^{R}_{2}$ &
$q$ & $\overline{t}_{2}y$ &
$q$ & $t^{R}_{2}$ &
$q$ & $t^{R}_{2}\vphantom{H^{H^{L^{L}}}}$ \\
 \hline
 39887 & 1173  &  39901 & 1181  &  39929 & 1175  &  39937 & 1165  &  39953 & 1173  &  39971 & 1178  \\
 39979 & 1169  &  39983 & 1180  &  39989 & 1169  &  39989 & 1172  &  40009 & 1173  &  40013 & 1179  \\
 40031 & 1179  &  40037 & 1167  &  40039 & 1166  &  40063 & 1177  &  40087 & 1164  &  40093 & 1174  \\
 40099 & 1174  &  40111 & 1170  &  40123 & 1169  &  40127 & 1176  &  40129 & 1174  &  40151 & 1172  \\
 40153 & 1185  &  40163 & 1174  &  40169 & 1172  &  40177 & 1166  &  40189 & 1173  &  40193 & 1172  \\
 40213 & 1173  &  40231 & 1173  &  40237 & 1173  &  40241 & 1181  &  40253 & 1173  &  40277 & 1176  \\
 40283 & 1181  &  40289 & 1185  &  40343 & 1175  &  40351 & 1180  &  40357 & 1184  &  40361 & 1178  \\
\hline
 \end{tabular}

 \newpage

\textbf{Table 1} Continue 21 \medskip

 \renewcommand{\arraystretch}{1.00}
 \begin{tabular}
 {@{}c@{\,\,\,}c|c@{\,\,\,}c|c@{\,\,\,}c|c@{\,\,\,}c|c@{\,\,\,}c|c@{\,\,\,}c@{}}
  \hline
$q$ & $t^{R}_{2}$ &
$q$ & $t^{R}_{2}$ &
$q$ & $t^{R}_{2}$ &
$q$ & $t^{R}_{2}$ &
$q$ & $\overline{t}_{2}y$ &
$q$ & $t^{R}_{2}\vphantom{H^{H^{L^{L}}}}$ \\
 \hline
 42131 & 1206  &  42139 & 1214  &  42157 & 1198  &  42169 & 1214  &  42179 & 1206  &  42181 & 1213  \\
 42187 & 1208  &  42193 & 1205  &  42197 & 1204  &  42209 & 1205  &  42221 & 1205  &  42223 & 1205  \\
 42227 & 1210  &  42239 & 1209  &  42257 & 1215  &  42281 & 1212  &  42283 & 1216  &  42293 & 1198  \\
 42299 & 1214  &  42307 & 1210  &  42323 & 1205  &  42331 & 1209  &  42337 & 1204  &  42349 & 1205  \\
 42359 & 1207  &  42373 & 1214  &  42379 & 1218  &  42391 & 1211  &  42397 & 1214  &  42403 & 1207  \\
 42407 & 1213  &  42409 & 1205  &  42433 & 1207  &  42437 & 1212  &  42443 & 1207  &  42451 & 1206  \\
 42457 & 1205  &  42461 & 1215  &  42463 & 1217  &  42467 & 1214  &  42473 & 1210  &  42487 & 1213  \\
\hline
 \end{tabular}

 \newpage

\textbf{Table 1} Continue 22 \medskip

 \renewcommand{\arraystretch}{1.00}
 \begin{tabular}
 {@{}c@{\,\,\,}c|c@{\,\,\,}c|c@{\,\,\,}c|c@{\,\,\,}c|c@{\,\,\,}c|c@{\,\,\,}c@{}}
  \hline
$q$ & $t^{R}_{2}$ &
$q$ & $t^{R}_{2}$ &
$q$ & $t^{R}_{2}$ &
$q$ & $\overline{t}_{2}yy$ &
$q$ & $t^{R}_{2}$ &
$q$ & $t^{R}_{2}\vphantom{H^{H^{L^{L}}}}$ \\
 \hline
 44357 & 1248  &  44371 & 1239  &  44381 & 1246  &  44383 & 1238  &  44389 & 1240  &  44417 & 1239  \\
 44449 & 1241  &  44453 & 1241  &  44483 & 1239  &  44491 & 1246  &  44497 & 1245  &  44501 & 1235  \\
 44507 & 1239  &  44519 & 1245  &  44531 & 1243  &  44533 & 1246  &  44537 & 1245  &  44543 & 1241  \\
 44549 & 1247  &  44563 & 1243  &  44579 & 1242  &  44587 & 1245  &  44617 & 1240  &  44621 & 1243  \\
 44623 & 1247  &  44633 & 1240  &  44641 & 1243  &  44647 & 1247  &  44651 & 1253  &  44657 & 1239  \\
 44683 & 1241  &  44687 & 1249  &  44699 & 1247  &  44701 & 1244  &  44711 & 1245  &  44729 & 1246  \\
 44741 & 1249  &  44753 & 1255  &  44771 & 1241  &  44773 & 1245  &  44777 & 1252  &  44789 & 1246  \\
\hline
 \end{tabular}


\newpage
\begin{thebibliography}{99}

\bibitem{abaV83} Abatangelo, V.: \emph{A class of complete
    $[(q+8)/3]$-arcs of $ \mathrm{PG}(2,q)$, with $q=2^{h}$ and
    $h$ ($\geq 6$) even.} Ars Combin. \textbf{16}, 103--111
    (1983)

\bibitem{tesi13}   Ali, A.H.: \emph{Classification of arcs in
    Galois plane of order thirteen}. Ph.D. Thesis, University
    of Sussex (1993).

\bibitem{Ball-SmallArcs} Ball, S.: \emph{On small complete arcs
    in a finite plane.} Discrete Math. \textbf{174}, 29--34
    (1997)
\bibitem{BDFKMP-Bulg2013} Bartoli, D., Davydov, A.A., Faina,
    G.,  Kreschuk, A.A., Marcugini, S., Pambianco,~F.: \emph{Two types of
    upper bounds on the smallest size of a complete arc in the
    plane $PG(2,q),$} in Proc. VII Int.
    Workshop on Optimal Codes and Related Topics,
    OC2013, Albena, Bulgaria, 2013, 19--25;
    available online at http://www.moi.math.bas.bg/oc2013/a3.pdf

\bibitem{BDFKMP-ArXiv} Bartoli, D., Davydov, A.A., Faina, G.,
     Kreschuk, A.A., Marcugini, S., Pambianco,~F.: \emph{Tables of sizes of small complete arcs in
the plane $PG(2,q)$, $q\le360007$}, (2013); available online at
http://arxiv.org/abs/1312.2155

\bibitem{BDFMP2012} Bartoli, D., Davydov, A.A., Faina, G.,
    Marcugini, S., Pambianco, F.: \emph{On sizes of complete
    arcs in $PG(2,q).$} Discrete Math \textbf{312}, 680--698
    (2012)

\bibitem{BDFMP-JG2013} Bartoli, D., Davydov, A.A., Faina, G.,
    Marcugini, S., Pambianco, F.: \emph{New upper bounds on the
    smallest size of a complete arc in a finite Desarguesian
    projective plane.} J. Geom.  {\bf 104}, 11--43 (2013)

\bibitem{BDFMP2013}  Bartoli, D., Davydov, A.A., Faina, G.,
    Marcugini, S., Pambianco, F.: \emph{A $3$-cycle
    construction of complete arcs sharing $(q+3)/2$ points with
    a conic.} Adv. Math. Commun. {\bf 7}(3),  319--334 (2013)

 \bibitem{BDFMP-ArXivFOP}   Bartoli, D., Davydov, A.A., Faina,
    G., Marcugini, S., Pambianco, F.:  \emph{Tables of
sizes of small complete arcs in the plane $PG(2, q)$, $q\le
190027$, obtained by an algorithm with fixed order of points
(FOP),} 2014, arXiv:1404.0469 [math.CO],
http://arxiv.org/abs/1404.0469

\bibitem{BDFMP_JGsubmit}  Bartoli, D., Davydov, A.A., Faina,
    G., Marcugini, S., Pambianco, F.: \emph{New types of estimates for the smallest size
    of complete arcs in a finite Desarguesian projective plane.}
    J. Geom.,
2014, DOI 10.1007/s00022-014-0224-4

\bibitem{BDMP_ACCT2012}  Bartoli, D., Davydov, A.A.,
    Marcugini, S., Pambianco, F.: \emph{New type of estimate for the smallest size of
complete arcs in PG(2; q).} In: Proceedings XIII International
Workshop on Algebraic and Combinatorial Coding Theory,
ACCT2012, Pomorie, Bulgaria, 2012, 67-72.

\bibitem{BFMP2013b} Bartoli, D., Faina,  G., Marcugini, S.,
    Pambianco,  F.,  Davydov,  A.A.: \emph{A new algorithm and
    a new type of estimate for the smallest size of complete
    arcs in $PG(2,q)$}. Electron. Notes Discrete Math. {\bf
    40}, 27--31 (2013)

\bibitem{BFMP2013} Bartoli, D., Faina, G., Marcugini, S.,
    Pambianco, F.: \emph{On the minimum size of complete arcs
    and minimal saturating sets in projective planes}. J.
    Geom. \textbf{104}, 409--419 (2013)

\bibitem{Blokhuis}  Blokhuis A., \emph{Blocking sets in
    desarguesian
    planes,} Erd\H{o}s is eighty, Bolyai society mathematical studies,  eds
     Mikl\'{o}s D.,  S\'{o}s V.T.,  Sz\'{o}nyi T., \textbf{2},    133--155 (1996).

\bibitem{CohLitZemGreedy} Cohen, G.,  Litsyn, S., Zemor,
    G.:\emph{On greedy algorithms in coding theory}. IEEE
    Trans. Inform. Theory \textbf{42}, 2053--2057 (1996)

\bibitem{DFMP-JG2005} Davydov, A.A., Faina, G., Marcugini, S.,
    Pambianco, F.: \emph{Computer search in projective planes
    for the sizes of complete arcs.} J. Geom.  \textbf{82},
    50--62 (2005)

\bibitem{DFMP-JG2009} Davydov, A.A., Faina, G., Marcugini, S.,
    Pambianco, F.: \emph{On sizes of complete caps in
    projective spaces $\mathrm{PG}(n,q)$ and arcs in planes $
    \mathrm{PG}(2,q)$.} J. Geom.  \textbf{94},  31--58  (2009)

\bibitem{DGMP-JCD} Davydov, A.A., Giulietti, M.,  Marcugini,
    S., Pambianco, F.: \emph{New inductive constructions of
    complete caps in $\mathrm{PG}(N,q)$, $q$ even.} J. Comb.
    Des. \textbf{18}, 176--201 (2010)

\bibitem{DGMP-AMC} Davydov, A.A., Giulietti, M., Marcugini, S.,
    Pambianco, F.: \emph{Linear nonbinary covering codes and
    saturating sets in projective spaces.} Adv. Math. Commun.
    \textbf{5}, 119--147 (2011)

\bibitem{DMP-JG2004} Davydov, A.A., Marcugini, S., Pambianco,
    F.: \emph{Complete caps in projective spaces
    $\mathrm{\mathrm{PG}}(n,q)$.} J. Geom. \textbf{80}, 23--30
    (2004)

\bibitem{faina15}  Faina,  G., Marcugini, S., Milani, A.,
    Pambianco,  F.: \emph{The spectrum of the values $k$ for
    which there exists a complete $k$-arc in $PG(2,q)$ for q
    $\leq 23$}. Ars Combin. {\bf 47}, 3--11  (1997)

\bibitem{FP} Faina, G., Pambianco, F.: \emph{On the spectrum of
    the values $k$ for which a complete $k$-cap in
    $\mathrm{PG}(n,q)$ exists.} J. Geom. \textbf{62}, 84--98
    (1998)

\bibitem{FPDM} Faina, G., Pambianco, F.: \emph{On some 10-arcs
    for deriving the minimum order for complete arcs in small
    projective planes.} Discrete Math. \textbf{208--209},
    261--271 (1999)

\bibitem{GacsSzonyi} G\'{a}cs, A., Sz\H{o}nyi, T.: \emph{Random
    constructions and density results.} Des. Codes Cryptogr.
    \textbf{47}, 267--287 (2008)

\bibitem{MaPa-private} Giulietti, M., Korchm\'{a}ros, G.,
    Marcugini, S., Pambianco, F.: \emph{Transitive
    $\mathbf{A_{6}}$-invariant $k$-arcs in $PG(2,q)$.}
    http://arxiv.org/abs/1108.0358

\bibitem{GiulUghi} Giulietti, M., Ughi, E.: \emph{A small
    complete arc in $\mathrm{PG}(2,q),$ $q=p^{2},$ $p\equiv
    3\mbox{ (mod }4)$.} Discrete Math. \textbf{208--209},
    311--318 (1999)

\bibitem{faina24}   Gordon, C.E.: \emph{Orbits of arcs in
    $PG(N,K)$ under projectivities.} Geom. Dedicata {\bf 42},
    187--203  (1992)

\bibitem{Hartman-Haskin} Hartman, A., Raskin, L.:
    \emph{Problems and algorithms for covering arrays.}
    Discrete Math. \textbf{284}, 149--156 (2004)


\bibitem{HirsSurvey83} Hirschfeld, J.W.P.: \emph{Maximum sets
    in finite projective spaces.} In: Lloyd, E.K. (ed.) Surveys
    in Combinatorics, London Math. Soc. Lecture Note Ser.
    \textbf{82}, pp. 55--76. Cambridge University Press,
    Cambridge (1983)

\bibitem{HirsBook} Hirschfeld, J.W.P.: \emph{Projective
    geometries over finite fields, 2nd edn.} Clarendon Press,
    Oxford (1998)

\bibitem{6}   Hirschfeld, J.W.P.,  Sadeh,  A.: \emph{The
    projective plane over the field of eleven elements.} Mitt.
    Math. Sem. Giessen {\bf 164}, 245--257  (1984)

\bibitem{HirsStor-2001} Hirschfeld, J.W.P., Storme, L.:
    \emph{The packing problem in statistics, coding theory and
    finite geometry: update 2001.} In:Blokhuis,  A.,
    Hirschfeld, J.W.P., Jungnickel, D., Thas, J.A. (eds.)
    Finite Geometries, Developments of Mathematics, vol. 3,
    Proc. of the Fourth Isle of Thorns Conf., Chelwood Gate,
    2000, pp. 201--246. Kluwer Academic Publisher, Boston
    (2001)

\bibitem{Keri} Keri, G.: \emph{Types of superregular matrices
    and the number of $n$-arcs and complete $n$-arcs in
    $\mathrm{PG}(r,q)$.} J. Comb. Des. \textbf{14}, 363--390
    (2006)

\bibitem{KV} Kim, J.H., Vu, V.: \emph{Small complete arcs in
    projective planes.} Combinatorica \textbf{23}, 311--363
    (2003)

\bibitem{korG83a} Korchm\'{a}ros, G.: \emph{New examples of
    complete $k$-arcs in $ \mathrm{PG}(2,q)$.} Europ. J.
    Combin. \textbf{4}, 329--334 (1983)

\bibitem{faina40}   Lisonek, P.: \emph{Computer-assisted
    studies in algebraic combinatorics.} Ph.D. thesis, Research
    Institute for Symbolic Computation, J. Kepler Univ. Linz,
    1994, RISC-Linz Report Series No. 94-68.

\bibitem{LombRad} Lombardo-Radice, L.: \emph{Sul problema dei
    k-archi completi di $ S_{2,q}$.} Boll. Unione Mat. Ital.
    \textbf{11}, 178--181 (1956)

\bibitem{MMP-q29} Marcugini, S., Milani, A., Pambianco, F.:
    \emph{Minimal complete arcs in $\mathrm{PG}(2,q),$ $q\leq
    29$.} J. Combin. Math. Combin. Comput. \textbf{47}, 19--29
    (2003)

\bibitem{art25} Marcugini, S., Milani A.,  Pambianco,  F.:
    \emph{Complete arcs in $PG(2,25)$: The spectrum of the
    sizes and the classification of the smallest complete
    arcs.} Discrete Math. {\bf 307}, 739--747  (2007)


\bibitem{Ost} \"{O}sterg\aa rd, P.R.J.: \emph{Computer search
    for small complete caps.} J. Geom. \textbf{69}, 172--179
    (2000)

\bibitem{PBFM2013} Pambianco, F., Bartoli, D., Faina, G.,
    Marcugini, S.:  \emph{Classification of the smallest
    minimal $1$-saturating sets in $PG(2,q)$, $q\leq23$}.
    Electron. Notes Discrete Math. {\bf 40}, 229--233 (2013)

\bibitem{PDBGM2013} Pambianco, F., Davydov, A.A., Bartoli, D.,
    Giulietti, M.,  Marcugini, S.: \emph{A note on multiple
    coverings of the farthest-off points}. Electron. Notes
    Discrete Math.  {\bf 40}, 289--293 (2013)


\bibitem{pelG77} Pellegrino, G.: \emph{Un'osservazione sul
    problema dei $k$-archi completi in $S_{2,q}$, con $q\equiv
    1$ $(\mbox{mod} 4)$.} Atti Accad. Naz. Lincei Rend.
    \textbf{63},  33--44 (1977)

\bibitem{pelG92} Pellegrino, G.: \emph{Sugli archi completi dei
    piani $PG(2,q)$, con $q$ dispari, contenenti
    $(q+3)/2$ punti di una conica.} Rend. Mat. \textbf{12},
    649--674 (1992)

\bibitem{pel93} Pellegrino, G.: \emph{Archi completi,
    contenenti $(q+1)/2$ punti di una conica, nei piani di
    Galois di ordine dispari.} Rend. Circ. Mat. Palermo (2)
    \textbf{62}, 273--308 (1993)

\bibitem{Polv} Polverino, O.:  \emph{Small minimal blocking
    sets and complete $k$ -arcs in $\mathrm{PG}(2,p^{3})$.}
    Discrete Math. \textbf{208--209}, 469--476 (1999).

\bibitem{SegreLeGeom} Segre, B.: \emph{Le geometrie di Galois.}
    Ann. Mat. Pura Appl. \textbf{48}, 1--97 (1959)

\bibitem{segB62b} Segre, B.: \emph{Ovali e curve $\sigma $ nei
    piani di Galois di caratteristica due.} Atti Accad. Naz.
    Lincei Rend. \textbf{32}, 785--790 (1962)

\bibitem{SegreIntrodGalGeom} Segre, B.: \emph{Introduction to
    Galois geometries.} Atti Accad. Naz. Lincei Mem.
    \textbf{8}, 133--236 (1967)

\bibitem{SZ} Sz\H{o}nyi, T.: \emph{Small complete arcs in
    Galois planes.} Geom. Dedicata \textbf{18}, 161--172 (1985)

\bibitem{szoT87a} Sz\H{o}nyi, T.: \emph{Note on the order of
    magnitude of $k$ for complete $k$-arcs in
    $PG(2,q)$.} Discrete Math. \textbf{66}, 279--282
    (1987)

\bibitem{szoT89survey} Sz\H{o}nyi, T.: \emph{Complete arcs in
    Galois planes: survey.} Quaderni del Seminario di Geometrie
    Combinatorie, Universit\`{a} degli studi di Roma, La
    Sapienza, \textbf{94} (1989).

\bibitem{szoT93} Sz\H{o}nyi, T.:  \emph{Arcs, caps, codes and
    3-independent subsets.} In: Faina, G., Tallini, G. (eds.)
    Giornate di Geometrie Combinatorie, Universit\`{a} degli
    studi di Perugia, pp. 57--80. Perugia, (1993)

\bibitem{Szonyi1997surveyCurves} Sz\H{o}nyi, T.:  \emph{Some
    applications of algebraic curves in finite geometry and
    combinatorics.} In: Bailey, R.A. ed. Surveys in
    Combinatorics, pp. 198--236. Cambridge University Press,
    Cambridge (1997)

\bibitem{tallini}  Tallini, G.: \emph{Le geometrie di Galois e
    le loro applicazioni alla statistica e alla teoria delle
    informazioni.} Rend. Mat. Appl. {\bf 19}, 379--400 (1960)

\bibitem{thaJ92d} Thas, J.A.: \emph{M.D.S. codes and arcs in
    projective spaces: a survey.} Le Matematiche (Catania)
    \textbf{47}, 315--328 (1992)

\bibitem{Voloch90} Voloch, J.F.: \emph{On the completeness of
    certain plane arcs II}. European J. Combin. \textbf{11},
    491--496 (1990)

\end{thebibliography}
\end{document}